\documentclass[12pt,reqno]{amsart}
\usepackage{amscd,amsmath,amsopn,amssymb,amsthm,multicol}
\usepackage{hyperref}
\usepackage[all]{xy}
\usepackage{amscd}
\usepackage{graphicx}

\DeclareMathOperator{\Ad}{Ad}
\DeclareMathOperator{\ad}{ad}
\DeclareMathOperator{\Aut}{Aut}

\DeclareMathOperator{\Ric}{Ric}

\newcommand{\fr}{\mathfrak}
\newcommand{\al}{\alpha}
\newcommand{\be}{\beta}

\newtheorem{theorem}{Theorem}

\newtheorem{prop}{Proposition}

\begin{document}

\title
{Homogeneous Einstein metrics on the generalized flag manifold $Sp(n)/(U(p)\times U(n-p))$}

\author{Andreas Arvanitoyeorgos, Ioannis Chrysikos and Yusuke Sakane}
\address{University of Patras, Department of Mathematics, GR-26500 Rion, Greece}
\email{arvanito@math.upatras.gr}
\email{xrysikos@master.math.upatras.gr}
\address{Osaka University, Department of Pure and Applied Mathematics, Graduate School of Information Science and Technology,
Osaka 560-043, Japan}
\email{sakane@math.sci.osaka-u.ac.jp}
 
\medskip
\noindent
\thanks{The third author was supported by Grant-in-Aid
for Scientific Research (C) 21540080}

\begin{abstract}
 We find the precise number of non-K\"ahler $Sp(n)$-invariant Einstein metrics on the generalized flag manifold
$M=Sp(n)/(U(p)\times U(n-p))$ with $n\geq 3$ and $1\leq p\leq n-1$.
We use an analysis on parametric systems of polynomial equations and we give some insight towards the study of such systems. 

 \medskip
\noindent 2000 {\it Mathematics Subject Classification.} Primary 53C25; Secondary 53C30, 12D05, 65H10

\medskip
\noindent {\it Keywords}:  homogeneous  manifold, Einstein metric,  generalized flag manifold,  algebraic system
of equations, resultant, Gr\"obner basis.

\end{abstract}

\maketitle

\section*{Introduction}
\markboth{Andreas Arvanitoyeorgos, Ioannis Chrysikos and Yusuke Sakane}{Invariant Einstein metrics on the generalized flag manifold 
$Sp(n)/(U(p)\times U(n-p))$}

A Riemannian metric $g$    is called {\it Einstein} if the Ricci tensor
 $\Ric_{g}$ satisfies the equation  
 ${\rm  Ric}_{g}=e\cdot g$, 
  for some $e\in\mathbb{R}$.     When $M$ is compact, 
  Einstein metrics of   volume 1 
  can be characterized variationally  as the critical 
  points of the scalar curvature functional 
    $T(g)=\int_{M}S_{g}d {\rm vol}_{g}$ 
  on the space  $\mathcal{M}_{1}$ of Riemannian metrics of 
  volume 1. If  $M=G/K$  is a compact homogeneous space, a $G$-invariant
Einstein metric  is precisely a critical 
point  of $T$ restricted to
the set  of $G$-invariant metrics of volume 1.
As a consequence, the Einstein equation reduces to a system 
of non-linear algebraic equations, 
which is still very complicated but more manageable, 
and sometimes can be solved explicity. 
 Thus most known examples of Einstein manifolds are homogeneous.  
   
 A generalized 
  flag manifold is an adjoint orbit  of a compact semisimple Lie group $G$, 
  or equivalently a compact homogeneous space  of the form  $M=G/K=G/ C(S)$, 
   where $C(S)$ is the centralizer  of a torus $S$ in $G$.      
  Einstein metrics on generalized flag manifolds have been studied by several authors
  (Alekseevsky, Arvanitoyeorgos, Kimura, Sakane, Chrysikos, Negreiros).
  
  Eventhough the problem of finding all invariant 
Einstein metrics on $M$  can be facilitated by use of certain
theoretical results (e.g. the work \cite{Gr} on 
the total number of $G$-invariant complex Einstein metrics),
it still remains a difficult one, especially when the number 
of isotropy summands increases.  This difficulty also increases 
when we pass from exceptional flag manifolds to classical flag 
manifolds, because in the later case the Einstein equation reduces 
to a parametric system.

  In two recent works \cite{Chry3} and \cite{Chry4}
  all invariant 
Einstein metrics were found for all generalized flag manifolds with
 four isotropy summands,  but a partial answer was given for the space $Sp(n)/(U(p)\times U(n-p))$.

We summarize  the results obtained  about this space.

\begin{theorem}\label{The2}{\textnormal{(\cite{Chry3})}}
The flag manifold $Sp(n)/(U(p)\times U(n-p))$ ($n\geq 2$ and $1\leq p\leq n-1$)  
    admits at least four 
    $Sp(n)$-invariant Einstein metrics, which are K\"ahler. 
\end{theorem}

\begin{theorem}\label{The4}{\textnormal{(\cite{Chry3})}}
The flag manifold $Sp(2p)/(U(p)\times U(p))$ ($p\geq 1$)  
    admits precisely six 
    $Sp(2p)$-invariant Einstein metrics.  There are four isometric K\"ahler-Einstein 
    metrics, and two   non-K\"ahler  Einstein metrics. 
\end{theorem}

In the present paper we find all $Sp(n)$-invariant Einstein 
metrics on the flag manifold $Sp(n)/(U(p)\times U(n-p))$, by using
an approach similar to the one used in \cite{Chry4}.

The Einstein equation reduces to polynomial systems whose
coefficients involve parameters, so a demanding task is 
to show  existence and uniqueness of solutions of such systems. 
 
Our main result is the following:

\medskip
{ \sc{Main Theorem.}}
    {\it The generalized flag manifold $M=Sp(n)/(U(p)\times U(n-p))$ with $n\geq 3$ and $1\leq p\leq n-1$,
    admits precisely two non-K\"ahler $Sp(n)$-invariant Einstein metrics.}

  \section{The Einstein equation for  flag manifolds}

Let $M=G/K=G/C(S)$ be a generalized flag manifold of a compact  
 simple Lie group $G$, where   $K=C(S)$ is the centralizer 
   of a torus $S$ in $G$.      
  Let   $o=eK$  be the identity coset of $G/K$. 
We denote by $\fr{g}$ and $\fr{k}$ the corresponding Lie 
algrebras of $G$ and $K$. Let $B$ denote the  
Killing form of $\fr{g}$.  Since $G$ is compact and  simple,  $-B$ is a positive 
definite inner product on $\fr{g}$.  With repsect to $-B$ we consider the orthogonal 
decomposition $\fr{g}=\fr{k}\oplus\fr{m}$. 
This is a reductive decomposition 
of $\fr{g}$, that is $\Ad(K)\fr{m}\subset\fr{m}$, and as usual we identify the tangent 
space $T_{o}M$   with $\fr{m}$.  Since $K=C(S)$, the isotropy group $K$ is connected and   the relation  $\Ad(K)\fr{m}\subset\fr{m}$ is equivalent with  $[\fr{k}, \fr{m}]\subset\fr{m}$.  Thus,  for a flag manifold $M=G/K$ the notion of $\Ad(K)$-invariant and $\ad(\fr{k})$-invariant is equivalent.

Let  $\chi : K\to \Aut(T_{o}M)$ 
be the isotropy representation of $K$ on $T_{o}M$. 
Since $\chi$ is equivalent to the adjoint representation   
of $K$ restricted on $\fr{m}$,  the set of  all $G$-invariant symmetric covariant  2-tensors on $G/K$  
can be identified with the set of all $\Ad(K)$-invariant symmetric bilinear forms on $\fr{m}$.  
In particular, the set of  $G$-invariant metrics on $G/K$ 
is identified with the set of $\Ad(K)$-invariant inner products  on $\fr{m}$.  

Let $\fr{m}=\fr{m}_1\oplus\cdots\oplus\fr{m}_{s}$ be a $(-B)$-orthogonal 
$\Ad(K)$-invariant decomposition of $\fr{m}$ into  pairwise inequivalent  irreducible 
$\Ad(K)$-modules $\fr{m}_{i}$  $(i=1, \ldots, s)$.  Such a decomposition always exists and can be expressed in terms of $\fr{t}$-roots (cf. \cite{AP}, \cite{Chry3}).
Then, a $G$-invariant Riemannian metric  on $M$ (or equivalently, an $\Ad(K)$-invariant inner product $\langle \ , \ \rangle$ on $\fr{m}=T_{o}M$) is given by 
\begin{equation}\label{Inva}
g=\langle \ , \ \rangle =x_1\cdot (-B)|_{\fr{m}_1}+\cdots+x_s\cdot (-B)|_{\fr{m}_s},
\end{equation}
where  $(x_1, \ldots, x_s)\in\mathbb{R}^{s}_{+}$.  Since $\fr{m}_{i}\neq\fr{m}_{j}$ as $\Ad(K)$-representations, any $G$-invariant metric on $M$ has the above form. 

Similarly, the Ricci tensor $\Ric_{g}$ of a $G$-invariant metric $g$ on $M$, as a symmetric covariant 2-tensor on $G/K$ 
  is given by
  \[
\Ric_{g} = r_1 x_1\cdot (-B)|_{\fr{m}_1}+\cdots+r_s x_s\cdot (-B)|_{\fr{m}_s}, 
\]
where $r_{1}, \ldots, r_{s}$  are the components of the Ricci 
tensor on each $\fr{m}_{i}$, that is $\Ric_{g}|_{\fr{m}_{i}}=r_{i} x_i\cdot (-B)|_{\fr{m}_{i}}$.  These components have o useful description in terms of the structure constants $[ijk]$ first introduced in \cite{Wa2}.   
   Let $\{X_{\al}\}$ be a $(-B)$-orthonormal basis adapted to the 
   decomposition of $\fr{m}$, that is  $X_{\al}\in \fr{m}_{i}$ for some $i$, 
   and $\al<\be$ if $i<j$ (with $X_{\al}\in \fr{m}_{i}$ and $X_{\be}\in\fr{m}_{j}$).  
   Set $A_{\al\be}^{\gamma}=B([X_{\al}, X_{\be}], X_{\gamma})$ so 
   that $[X_{\al}, X_{\be}]_{\fr{m}}=\sum_{\gamma}A_{\al\be}^{\gamma}X_{\gamma}$, 
   and   $[ijk]=\sum(A_{\al\be}^{\gamma})^{2}$, where the sum is taken over all 
   indices $\al, \be, \gamma$ with $X_{\al}\in \fr{m}_{i}, X_{\be}\in\fr{m}_{j}, X_{\gamma}\in\fr{m}_{k}$
   (where $[ \ , \ ]_{\fr{m}}$ denotes the $\fr{m}$-component).  Then $[ijk]$ is nonnegative, symmetric in all 
   three entries, and independent of the $(-B)$-orthonormal bases choosen for 
   $\fr{m}_{i}, \fr{m}_{j}$ and $\fr{m}_{k}$ (but it depends on the choise 
   of the decomposition of $\fr{m}$).

 \begin{prop}\label{Ricc}{\textnormal{(\cite{SP})}} 
Let $M=G/K$ be a generalized flag manifold of a compact simple Lie group $G$ and let
 $\fr{m}=\bigoplus_{i=1}^{s}\fr{m}_{i}$ be a decomposition of $\fr{m}$ into pairwise
 inequivalent irreducible $\Ad(K)$-submodules.  
Then the components $r_{1}, \ldots, r_{s}$ of the Ricci tensor of a $G$-invariant 
metric  \em{(\ref{Inva})} on $M$ are given by
   \[
   r_{k}=\frac{1}{2x_{k}}+\frac{1}{4d_{k}}\sum_{i, j}\frac{x_{k}}{x_{i}x_{j}}[ijk]-\frac{1}{2d_{k}}\sum_{i, j}\frac{x_{j}}{x_{k}x_{i}}[kij], \qquad (k=1, \ldots, s).
\]
\end{prop}

 In wiew of Proposition \ref{Ricc}, a $G$-invariant metric $g=(x_1, \ldots, x_{s})\in\mathbb{R}^{s}_{+}$ on $M$, is an Einstein metric with Einstein constant $e$, if and only if it is a positive real solution of the system 
 \[
 \frac{1}{2x_{k}}+\frac{1}{4d_{k}}\sum_{i, j}\frac{x_{k}}{x_{i}x_{j}}[ijk]-\frac{1}{2d_{k}}\sum_{i, j}\frac{x_{j}}{x_{k}x_{i}}[kij]=e, \quad 1\leq k\leq s.
  \]

\section{The generalized flag manifold $Sp(n)/(U(p)\times U(n-p))$}

  We review some results related to the generalized flag manifold $M=G/K=Sp(n)/(U(p)\times U(n-p))$ ($n\ge 3, \ 1\le p\le n-1$)
  obtained in \cite{Chry3}.
  Its corresponding painted Dynkin diagram is given by

\hspace{4.5cm}
\begin{picture}(160,45)(-15,-25)
\put(0, 0){\circle{4}}
\put(0,8.5){\makebox(0,0){$\al_1$}}
\put(0,-8){\makebox(0,0){2}}
\put(2, 0.3){\line(1,0){14}}
\put(18, 0){\circle{4}}
\put(18,8.5){\makebox(0,0){$\al_2$}}
\put(18,-8){\makebox(0,0){2}}
\put(20, 0.3){\line(1,0){10}}
\put(40, 0){\makebox(0,0){$\ldots$}}
\put(50, 0.3){\line(1,0){10}}
\put(60, 0){\circle*{4.4}}
\put(60,8.5){\makebox(0,0){$\al_p$}}
\put(60,-8){\makebox(0,0){2}}
\put(60, 0.3){\line(1,0){10}}
\put(80, 0){\makebox(0,0){$\ldots$}}
\put(90, 0.3){\line(1,0){10}}
\put(102, 0){\circle{4}}
\put(102,8.5){\makebox(0,0){$\al_{n-1}$}}
\put(102,-8){\makebox(0,0){2}}
\put(107.2, 1.2){\line(1,0){14.8}}
\put(107.2, -1){\line(1,0){14.8}}
\put(103.46, -2){\scriptsize $<$}
\put(123.5, 0){\circle*{4}}
\put(124,8.5){\makebox(0,0){$\al_n$}}
\put(123,-8){\makebox(0,0){1}}
\end{picture}

 \medskip

  The isotropy representation of $M$ decomposes into a direct sum $\chi = \chi _1\oplus\chi _2\oplus\chi _3\oplus\chi _4$, which
  gives rise to a decomposition $\fr{m}=\fr{m}_1\oplus\fr{m}_2\oplus\fr{m}_3\oplus\fr{m}_4$ of $\fr{m}=T_{o}M$
  into   four irreducible inequivalent  $\ad(\fr{k})$-submodules. 
  The dimensions $d_i =\dim \fr{m}_i\ (i=1,2,3,4)$ of these
  submodules can be obtained by use of Weyl's formula \cite[pp.~204-205, p.~210]{Chry3}    
  and are given by
  $$d_1=2p(n-p),  \ d_2=(n-p)(n-p+1), \  d_3=2p(n-p),\   d_4=p(p+1).
  $$

According to (\ref{Inva}), a $G$-invariant metric on $M=G/K$ is given by
\begin{equation}\label{metrI}
\left\langle \ , \ \right\rangle=x_1\cdot (-B)|_{\fr{m}_1}+
x_2\cdot (-B)|_{\fr{m}_2}+x_3\cdot (-B)|_{\fr{m}_3}+x_4\cdot (-B)|_{\fr{m}_4},
\end{equation}
 for positive real numbers $x_1, x_2, x_3, x_4$.
We will denote such metrics by $g=(x_1, x_2, x_3, x_4)$.

 It is known (\cite{N}) that if $n\ne 2p$ then $M$ admits two non-equivalent $G$-invariant complex structures $J_1, J_2$, and thus two non-isometric  K\"ahler-Einstein metrics which are given (up to scale) by (see also \cite[Theorem 3]{Chry3}) 
 \begin{equation}\label{KE}
 \begin{tabular}{l}  
 $g_1= (n/2, \ n+p+1, \ n/2+p+1, \ p+1)$ \\
 $g_2= (n/2, \ n-p+1, \ 3n/2-p+1, \ 2n-p+1)$.
 \end{tabular}
 \end{equation}
 If $n=2p$ then $M$ admits a unique $G$-invariant complex structure with
 corresponding K\"ahler-Einstein metric (up to scale) 
  $g= (p, \ p+1, \ 2p+1, \ 3p+1)$ (cf. also \cite[Theorem 10]{Chry3} where all isometric 
  K\"ahler-Einstein metrics are listed).

 The Ricci tensor of $M$ is given as follows:
 
 \begin{prop}\label{componentsII} (\cite{Chry3})  The  components $r_i$ of the Ricci tensor
  for a   $G$-invariant Riemannian metric on $M$ determined by {\rm (\ref{metrI})} are given as follows:

\begin{equation}\label{compIII}
 \ \ \left.
\begin{tabular}{l}
$r_1=\displaystyle\frac{1}{2x_1}+  \frac{c_{12}^3}{2d_1}\Big( \frac{x_1}{x_2x_3}- \frac{x_2}{x_1x_3}- \frac{x_3}{x_1x_2}\Big)+  \frac{c_{13}^4}{2d_1}\Big( \frac{x_1}{x_3x_4}- \frac{x_4}{x_1x_3}- \frac{x_3}{x_1x_4}\Big)$ \\
$r_2=\displaystyle\frac{1}{2x_2}+  \frac{c_{12}^3}{2d_2}\Big( \frac{x_2}{x_1x_3}- \frac{x_1}{x_2x_3}- \frac{x_3}{x_1x_2}\Big)$ \\
$r_3=\displaystyle\frac{1}{2x_3}+  \frac{c_{12}^3}{2d_3}\Big( \frac{x_3}{x_1x_2}- \frac{x_2}{x_1x_3}- \frac{x_1}{x_2x_3}\Big)+  \frac{c_{13}^4}{2d_3}\Big( \frac{x_3}{x_1x_4}- \frac{x_4}{x_1x_3}- \frac{x_1}{x_3x_4}\Big)$ \\
$r_4 = \displaystyle\frac{1}{2x_4}+  \frac{c_{13}^4}{2d_4}\Big( \frac{x_4}{x_1x_3}- \frac{x_3}{x_1x_4}- \frac{x_1}{x_3x_4}\Big),$ \\
\end{tabular}\right\}
\end{equation}
where $c_{12}^3=[123]$ and $c_{13}^4=[134]$.
\end{prop}

 By taking into account the explicit form of the K\"ahler-Einstein metrics above,  
 and substituting these in  (\ref{compIII}), 
 we find that the   values of the unknown triples $[ijk]$ are 
$\displaystyle{ c_{12}^3=\frac{p(n-p)(n-p+1)}{2(n+1)}}$  and $\displaystyle{ c_{13}^4=\frac{p(p+1)(n-p)}{2(n+1)}}$.

\medskip

A $G$-invariant metric $g=(x_1, x_2, x_3, x_4)$ on $M=G/K$ is Einstein if and only if, there is a positive
constant $e$ such that $r_1=r_2=r_3=r_4=e$, or equivalently
\begin{equation}\label{systemI}
   r_1-r_3=0, \quad r_1-r_2=0, \quad r_3-r_4=0. 
\end{equation} 

By substituting the values of $d_i\ (i=1,2,3,4)$ and $c_{12}^3, c_{13}^4$ into the components of the Ricci tensor, then
System (\ref{systemI}) is equivalent to the following equations:  

\medskip

\begin{equation}\label{syst5}
 \left.
 \begin{tabular}{r}
$(x_1-x_3)(x_1 x_2 + p x_1 x_2 + x_2 x_3 + p x_2 x_3 + x_1 x_4   + n x_1 x_4$\\
 $-  p x_1 x_4 - 2 x_2 x_4    - 2 n x_2 x_4 + x_3 x_4 + n x_3 x_4 - p x_3 x_4)=0$\\
$4(n+1) x_3 x_4 (x_2-x_1)+ (n+p+1)x_4(x_1^2-x_2^2)- (n-3p+1)x_3^2x_4$\\
$ +(p+1)x_2(x_1^2-x_3^2-x_4^2)=0$\\
$4(n+1)x_1x_2(x_4-x_3)+ (2n-p+1)x_2(x_3^2-x_4^2)+ (2n-3p-1)x_1^2x_2$\\
$ +(n-p+1)x_4(x_3^2-x_1^2-x_2^2)=0$\\
\end{tabular}\right\} 
  \end{equation}

\section{Proof of the Main Theorem}

  Consider   the equation $r_1- r_3 = 0$.  Then we have 
  \begin{eqnarray*}  & (x_1 - x_3) (x_1 x_2 + p x_1 x_2 + x_2 x_3 + p x_2 x_3 + x_1 x_4   + n x_1 x_4  & \\   & -  p x_1 x_4 - 2 x_2 x_4    - 2 n x_2 x_4 + x_3 x_4 + n x_3 x_4 - p x_3 x_4) =0.  &
 \end{eqnarray*}
We claim that, for  $x_1= x_3 $,  there are no Einstein metrics, and that, for the other case, there exist  Einstein metrics on $Sp(n)/( U(p) \times U( n -p))$.  

\medskip
  
{\bf Case 1.}  For  $x_1= x_3 $,  we put $x_1 =1$ and 
  we  get  the following system of  equations
  \begin{eqnarray} & & 
   (n+p+1){x_2}^2  +4 (n-p+1)-4 (n+1) {x_2}+(p+1) {x_2} {x_4} =0,  \label{1}  \\
 & &    (n-p+1)  {x_2} {x_4} + (2 n-p+1){x_4}^2 - 4 (n+1) {x_4} + 4 (p+1) = 0. \label{2}
\end{eqnarray} 

From equation (\ref{2}), we have 
 \begin{eqnarray} & & 
  x_2 =  \frac{{ -(2 n-p+1) x_4}^2 +4 (n+1) {x_4}-4
   (p+1)}{{x_4} (n-p+1)}. \label{3}
   \end{eqnarray} 
 
  Now we substitute equation  (\ref{3})  into the equation (\ref{1}), and we obtain the following equation:  
 \begin{eqnarray}
& & f_{n, p}(x_4) =   n (n+1)(2 n-p+1){x_4}^4-4 (n+1) \left(n^2+2 n
   p+n-p^2+p\right) {x_4}^3   \nonumber  \\ 
 &  &  +2 \left(n^3+9 n^2 p+7 n^2+4 n p^2+16 n p+8
   n-2 p^3+2 p^2+6 p+2\right){x_4}^2   \nonumber \\
 & & -  8 (n+1) (p+1)
    (n+3 p+1){x_4} +8 (p+1)^2 (n+p+1) = 0.\nonumber  
    \end{eqnarray}
   
   From  equation (\ref{2}), we have 
 \begin{eqnarray} & & 
  x_4= \frac{ -(n+p+1){x_2}^2+4(n+1){x_2}-4 (n-p+1)}{(p+1){x_2}}. \label{4}
   \end{eqnarray} 
   
    Now we substitute equation  (\ref{4})  into equation (\ref{1}), and we obtain the following equation:  
 \begin{eqnarray}
& & g_{n, p}(x_2) = n (n+1)(n+p+1){x_2}^4 -4 (n+1)  \left(2 n^2+2
   n-p^2-p\right) {x_2}^3   \nonumber   
\\
& &   +2 \left(12
   n^3-11 n^2 p+25 n^2-2 n p^2-20 n p+14
   n+2 p^3+2 p^2-6 p+2\right)  {x_2}^2  \nonumber   
 \\ 
 & &  - 8 (n+1) {x_2} (4
   n-3 p+1) (n-p+1)+8 (n-p+1)^2 (2 n-p+1)  = 0.  \nonumber   
    \end{eqnarray}

 We claim that, for $n/2 \leq p \leq n-1$,  
there are no positive solutions of the equation $f_{n, p}(x_4) = 0$. 
Note that $  g_{n, p}(x_4) = f_{n, n-p}(x_4)$, thus we also see that, for $1 \leq p \leq n/2$,  
there are no positive solutions of the equation $g_{n, p}(x_2) = 0 $.

    It is   
 \begin{eqnarray*}
 & & \frac{d f_{n, p}}{d x_4 }(x_4) =  4 n (n+1)(2 n-p+1) {x_4}^3
 -12 (n+1)\left(n^2+2 n p+n-p^2+p\right){x_4}^2\\
   & &  + 4 \left(n^3+9
   n^2 p+7 n^2+4 n p^2+16 n p+8 n-2 p^3+2
   p^2+6 p+2\right) {x_4} \\
   & & - 8 (n+1) (p+1) (n+3 p+1).
        \end{eqnarray*} 
   Note that the coefficient of ${x_4}^3$ is $4 n (n+1)(2 n-p+1) > 0$.

By evaluating $\displaystyle \frac{d f_{n, p}}{d x_4 }(x_4)$ at $\displaystyle x_4=  \frac{2(p-1)}{n}$, we have 
    \begin{eqnarray*}
 & & \frac{d  f_{n, p}}{d x_4 }\left( \frac{2(p-1)}{n} \right) = 
-\frac{8 (n-p+1) \left(2 n^3-3 n^2 p+13
   n^2-8 n p+16 n+2 p^3-6 p+4\right)}{n^2}.   \end{eqnarray*} 
Since we can write 
  \begin{eqnarray*}   & &  \ \ \   2 n^3-3 n^2 p+13
   n^2-8 n p+16 n+2 p^3-6 p+4 \\
  & &  = 
2 (n-p)^3+(3 p+13) (n-p)^2+(18 p+16)
   (n-p)+p^3+5 p^2+10 p+4, 
    \end{eqnarray*} 
   we see  that $\displaystyle  \frac{d  f_{n, p}}{d x_4 }\left( \frac{2(p-1)}{n} \right) < 0$. 

  By evaluating $\displaystyle \frac{d  f_{n, p}}{d x_4 }(x_4)$ at $\displaystyle x_4=  \frac{2(p+1)}{n}$, we see that 
    \begin{eqnarray*}
 & & \frac{d  f_{n, p}}{d x_4 }\left( \frac{2(p+1)}{n} \right) = 
-\frac{8(p+1) (n-p+1) \left(n^2-4 n p-4 n+2 p^2-2 p-4\right)}{n^2}.  
   \end{eqnarray*} 
For $\displaystyle \frac{n}{2} \leq p \leq n -1$,  we see that   
  \begin{eqnarray*}  & &  \ \ \  n^2-4 n p-4 n+2 p^2-2 p-4 
  = -\frac{n^2}{2}+2
   \left(p-\frac{n}{2}\right)^2+(-2 n-2)  \left(p-\frac{n}{2}\right)-5 n-4\\    
    &   & \leq  -\frac{n^2}{2} +2 \left( n -1 -\frac{n}{2} \right)^2 +(-2 n-2) \left(p-\frac{n}{2}\right) -5 n-4 \\    
    &   & = (-2 n-2) \left(p-\frac{n}{2}\right)-7 n-2 < 0,         \end{eqnarray*} 
and thus  we have  $\displaystyle  \frac{d f_{n, p}}{d x_4 }\left(\frac{2(p+1)}{n}  \right) > 0$. 
   Hence, for $\displaystyle \frac{n}{2} \leq p \leq n -1$,   
   the equation    $\displaystyle \frac{d f_{n, p}}{d x_4 }(x_4) = 0$ has  a real solution $u_1$   with 
   $$\displaystyle   \frac{ 2(p - 1)}{2 n} < u_1 < \frac{2(p + 1)}{n}. $$ 
   
 We claim that  the polynomial $ \displaystyle  \frac{d f_{n, p}}{d x_4 }$ of degree 3 is monotone increasing, and hence $f_{n, p}$ attains a local minimum at $x_4 = u_1$.    

Now the second derivative of $ f_{n, p}$ is given by 
  \begin{eqnarray*} 
 &   &  \ \  \ \frac{d^2 f_{n, p}}{d {x_4}^2 }(x_4) = 12 n (n+1)(2 n-p+1){x_4}^2 -24 (n+1)\left(n^2+2 n
   p+n-p^2+p\right){x_4} \\
  &   &  4 \left(n^3+n^2 (9
   p+7)+4 n \left(p^2+4 p+2\right)-2 (p+1)
   \left(p^2-2 p-1\right)\right). 
  \end{eqnarray*} 
     To see that  the second derivative of $ f_{n, p}$ is positive for $\displaystyle \frac{n}{2} \leq p \leq n -1$,  
 we note that the discriminant of the polynomial $\displaystyle \frac{d^2 f_{n, p}}{d {x_4}^2 }$ of degree 2 is given by 
      \begin{eqnarray*} 
        &   & 
 \left(12 (n+1) \left(n^2+2 n
   p+n-p^2+p\right)\right)^2 -12 n (n+1) (2 n-p+1)\times  
    \\
     &   & 
  4 \left(n^3+n^2 (9
   p+7)+4 n \left(p^2+4 p+2\right)-2 (p+1)
   \left(p^2-2 p-1\right)\right) \\
   &   & =  48 (n+1) \left(n^5-5 n^4 p-6 n^4+7 n^3
   p^2-4 n^3 p-14 n^3-4 n^2 p^3+20 n^2
   p^2+4 n^2 p \right.  \\
     &   & \left.  -9 n^2+n p^4-14 n p^3+13 n
   p^2+2 n p-2 n+3 p^4-6 p^3+3 p^2\right).  
    \end{eqnarray*}   
    
 We put 
      \begin{eqnarray*} &   &  
  h_{n, p} =  \left(n^5-5 n^4 p-6 n^4+7 n^3
   p^2-4 n^3 p-14 n^3-4 n^2 p^3+20 n^2
   p^2+4 n^2 p \right.  \\
     &   & \left.  -9 n^2+n p^4-14 n p^3+13 n
   p^2+2 n p-2 n+3 p^4-6 p^3+3 p^2\right).  
         \end{eqnarray*}   
  We consider $h_{n, p}$ as a polynomial of $p$ and we show that $ h_{n, p} < 0$ for $\displaystyle \frac{n}{2} \leq p \leq n -1$. 
 We have
  \begin{eqnarray*} &   &   h_{n, p }(p)= (n+3) p^4-2
   (n+3) (2 n+1) p^3
+ \left(7 n^3+20 n^2+13 n+3\right) p^2 \\
 &   & -n  \left(5 n^3+4 n^2-4 n-2\right) p+n (n+1)
   \left(n^3-7 n^2-7 n-2\right),   
       \end{eqnarray*}   
  \begin{eqnarray*} &   &    \frac{d h_{n, p }}{d {p} }(p)= 4 (n+3)
   p^3-6 (n+3) (2 n+1) p^2 \\
  &   &  + 2 \left(7 n^3+20 n^2+13 n+3\right) p-n
   \left(5 n^3+4 n^2-4 n-2\right)  
        \end{eqnarray*}     
   
   and     
        \begin{eqnarray*} 
        &   & \frac{d^2 h_{n, p }}{d p^2 }(p) = 
 12 (n+3) p^2-12 (n+3) (2 n+1) p +  2 \left(7 n^3+20 n^2+13 n+3\right) \\
   &   &  = 12 (n+3) (-n+p+1)^2-36 (n+3) (-n+p+1) + 2 \left(n^3-4 n^2+7 n+39\right) \\
    &   &   = 12 (n+3) (-n+p+1)^2-36 (n+3) (-n+p+1) \\
       &   &  + 2 (n-2)^3+4 (n-2)^2+6 (n-2)+90  > 0\end{eqnarray*}         
       
    for  $n \geq 2$  and  $   p \leq n-1$.      
    
    Thus $\displaystyle  \frac{d h_{n, p }}{d {p} } $ is  a monotone increasing function,  and we see that  
    $$\frac{d h_{n, p }}{d {p} }\left(\frac{2}{3} n \right)= \frac{1}{27} n \left(5 n^3+204 n^2+360
   n+162\right) > 0,  $$
    $$\frac{d h_{n, p }}{d {p} }\left(\frac{1}{2} n \right)= -\frac{1}{16} n \left(3 n^4+73 n^3+152
   n^2+116 n+32\right) <  0.  $$
   Thus the equation $\displaystyle  \frac{d h_{n, p }}{d {p} } = 0$ has a unique solution  $ \alpha$ with $ \displaystyle  \frac{1}{2}n  < \alpha  < \frac{2}{3} n$ and 
the function $h_{n, p}$  attains the minimum  only at $p =  \alpha$. 

Note that  $\displaystyle  h_{n, p } \left( \frac{n}{2}\right)  = -\frac{1}{16} n \left(3 n^4+73 n^3+152
   n^2+116 n+32\right) <  0$ and $\displaystyle   h_{n, p }(n-1)  = -2 \left(n^4+4 n^3+10 n^2+6 n-6\right)  <  0$. Thus we get that 
 $ h_{n, p} < 0$ for $\displaystyle \frac{n}{2} \leq p \leq n -1$.     
 
Since $f_{n, p }(0) = 8 (p+1)^2 (n+p+1) > 0$, in order to show that $ f_{n, p }(x_4) > 0 $ for $x_4 >0$,  we need to prove that
the local minimum $f_{n, p }(u_1)$  is positive.  

We consider the tangent lines $l_1, l_2, l_3$  of the curve $ f_{n, p }(x_4)$ at the points
$P_1$ with x-coordinate  $x_4 = 2 (p-1)/n$,  $P_2$ with $x_4 =2 (p + 1)/n$, 
and $P_3$ with  $x_4 =2 p /n$.  The equation of the line $l_1$ is given by 
 \begin{eqnarray*}
 & &  \ \ \  l_1(t) = 
\frac{16}{n^3} (n-p+1)^2 \left(2 n^2+4
   n+p^3-p^2-p+1\right) \\
& &   -\frac{8}{n^2}
   (n-p+1) \left(2 n^3-3 n^2 p+13 n^2-8 n
   p+16 n+2 p^3-6 p+4\right)
   \left(t-\frac{2 (p-1)}{n}\right), 
 \end{eqnarray*}  
 the equation of the line $l_2$ is given by 
   \begin{eqnarray*}
   & & \ \ \  l_2(t) = \frac{16}{n^3} (p+1)^3 (n-p+1)^2 \\
  & &    -\frac{8}{n^2}
   (p+1) (n-p+1) \left(n^2-4 n p-4 n+2
   p^2-2 p-4\right) \left(t-\frac{2
   (p+1)}{n}\right), \end{eqnarray*}  
 and    the equation of the line $l_3$ is given by 
   \begin{eqnarray*}
   & & \ \ \  l_3(t) = \frac{8}{n^3} \left(n^4-n^3 p+n^3+2 n^2 p^3-2 n^2
   p-4 n p^4+2 n p^3+2 n p^2+2 p^5-2
   p^4\right) \\  
    & &  -\frac{8}{n^2} (n-p+1)
   \left(n^3-n^2 p+n^2-2 n p^2-2 n p+2
   p^3\right) \left(t-\frac{2
   p}{n}\right). 
    \end{eqnarray*}

 \begin{minipage}{.54\linewidth}
Let $P_{13}$ be the point at which 
the tangent lines $l_1, l_3$  intersect and   $P_{23}$ the  point at which 
the tangent lines $l_2, l_3$  intersect. The coordinates $ ( \alpha_1, \beta_1 )$ of the point $P_{13}$ are  given by
Let $P_{13}$ be the point at which 
the tangent lines $l_1, l_3$  intersect and   $P_{23}$ the  point at which 
the tangent lines $l_2, l_3$  intersect. The coordinates $ ( \alpha_1, \beta_1 )$ of the point $P_{13}$ are  given by 
\end{minipage}
\begin{minipage}{.44\linewidth}
\begin{center}
\includegraphics[width=0.9\linewidth]{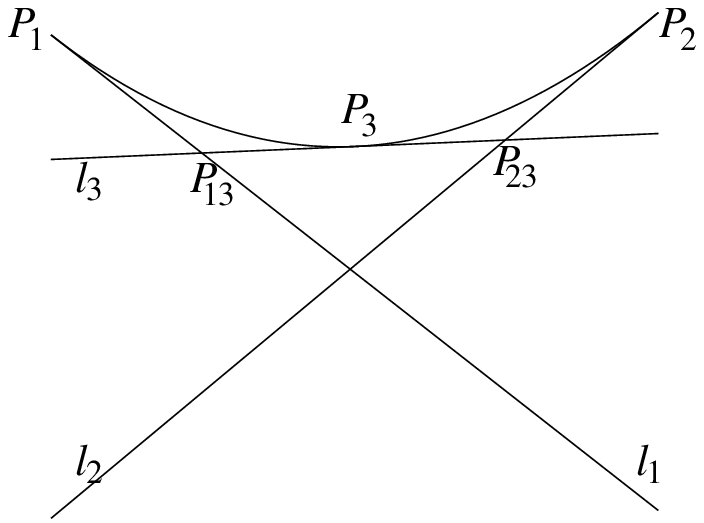} 
\end{center}
\end{minipage}
 \begin{eqnarray*}
 & & \alpha_1 =  \left(\left(2 p^2+3 p-1\right) (n-p)^3+\left(2
   p^3+15 p^2+13 p+4\right) (n-p)^2\right. \\
    & & +  \left(2
   p^4+15 p^3+39 p^2+36 p+12\right)
   (n-p)+(2 p+1) (n-p)^4 +2 p^4+13 p^3 \\
  & &  \left. +28 p^2+22 p+6\right)/(n (n-p+1) \left((n-p)^3+p (n-p)^2+\left(p^2+6 p+4\right) (n-p) \right. \\
 & & \left.  +(p+2) \left(p^2+4 p+2\right)\right)), \\
& &    \beta_1 = \left(8 (p+1) \left(n^6+2 n^5 p^2-3 n^5
   p+n^5-8 n^4 p^3+6 n^4 p^2-2 n^4 p-8
   n^4+14 n^3 p^4 \right. \right.  \\
 & &      +4 n^3 p^3+18 n^3 p^2+14
   n^3 p-14 n^3-12 n^2 p^5-16 n^2 p^4-4 n^2
   p^3+28 n^2 p^2+26 n^2 p \\
  & &  \left.   -6 n^2+4 n p^6+8
   n p^5-24 n p^4-28 n p^3+12 n p^2+12 n
   p+12 p^5-12 p^3\right)/ \\
& &     \left(n^3 \left((n-p)^3+p (n-p)^2+\left(p^2+6 p+4\right) (n-p)+(p+2) \left(p^2+4 p+2\right)\right)\right) \\
& & =  \left(8 (p+1) \left((n-p)^6+\left(2 p^2+3 p+1\right) (n-p)^5+\left(2
   p^3+6 p^2+3 p-8\right) (n-p)^4 \right. \right.  \\
   & & 
  +\left(2
   p^4+18 p^3+20 p^2-18 p-14\right)
   (n-p)^3+\left(2 p^5+17 p^4+48 p^3+22
   p^2-16 p-6\right)\times \\
   & &  \left.\left.  (n-p)^2+\left(3 p^5+19
   p^4+38 p^3+22 p^2\right)(n-p)+p^2 \left(p^3+6 p^2+12
   p+6\right)\right) \right)/ \\
& &     \left(n^3 \left((n-p)^3+p (n-p)^2+\left(p^2+6 p+4\right) (n-p)+(p+2) \left(p^2+4 p+2\right)\right)\right)
       \end{eqnarray*}  
and the coordinates $ ( \alpha_2, \beta_2 )$ of the point $P_{23}$ are  given by 
 \begin{eqnarray*}
 & & \alpha_2 = \left((2 p-1)(n-p)^4+ \left(2 p^2+25 p-15\right)
   (n-p)^3 \right.  \\
& &   +\left(2 p^3+37 p^2+35
   p-36\right) (n-p)^2 +\left(2 p^4+13
   p^3+49 p^2+4 p-28\right) (n-p)  \\
 & &  \left. +2 p^4+11 p^3+12 p^2-6 p-6 \right)/ \left(n (n-p+1)
  \left( (n-p)^3+ (p+12) (n-p)^2 \right. \right.\\
& &  \left.  \left. +\left(p^2+18 p+16\right) (n-p)+(p+2) \left(p^2+4 p+2\right)\right) \right), \\ 
 & &    \beta_2 =  8 \left(2 n^7-7 n^6 p+29 n^6+2 n^5
   p^3+9 n^5 p^2-66 n^5 p+79 n^5-8 n^4
   p^4+26 n^4 p^3\right.   \\
  & &     +28 n^4 p^2-178 n^4 p+84
   n^4+14 n^3 p^5-78 n^3 p^4+106 n^3 p^3+72
   n^3 p^2-192 n^3 p+38 n^3 \\
    & & -12 n^2 p^6+68
   n^2 p^5-172 n^2 p^4+144 n^2 p^3+82 n^2
   p^2-84 n^2 p+6 n^2+4 n p^7-20 n p^6 \\
   & &   \left. +88 n
   p^5-148 n p^4+64 n p^3+24 n p^2-12 n
   p-12 p^6+36 p^5-36 p^4+12
   p^3  \right)/
   \\
   & &   \left(n^3 \left( (n-p)^3+(p+12) (n-p)^2+\left(p^2+18 p+16\right)
   (n-p)+(p+2)
   \left(p^2+4 p+2\right) \right)\right)    \\
     & & = 8 \left(
     2 (n-p)^7+(7 p+29)
   (n-p)^6+ \left(2 p^3+9
   p^2+108 p+79\right) (n-p)^5 \right. \\
    & & 
   +\left(2
   p^4+36 p^3+133 p^2+217 p+84\right)
   (n-p)^4   \\
    & & 
+ \left(2 p^5+46 p^4+138 p^3+150
   p^2+144 p+38\right) (n-p)^3   \\
    & & \left.
+\left(2
   p^6+17 p^5+89 p^4+82 p^3+10 p^2+30
   p+6\right) (n-p)^2  \right. \\
    & & \left.
+ 3 \left((p-1)^4+10 (p-1)^3+37 (p-1)^2+44
   (p-1)+6\right) p^2 (n-p)  \right. \\
    & & \left.
+\left((p-1)^4+9 (p-1)^3+23
   (p-1)^2+13 (p-1)-8\right) p^2
   \right)/
   \\
   & &   \left(n^3 \left( (n-p)^3+(p+12) (n-p)^2+\left(p^2+18 p+16\right)
   (n-p)+(p+2)
   \left(p^2+4 p+2\right) \right)\right). 
    \end{eqnarray*}  
    Note that $\beta_1, \beta_2$ are positive for $ 1 \leq p \leq n-1$. 

Since $  f_{n, p}(x_4)$ is concave up,   we see that the curve $ ( x_4, f_{n, p}(x_4) ) $ for
$\displaystyle (\frac{2(p-1)}{n} \leq x_4 \leq  \frac{2 p}{n}  )$ lies inside the triangle given by the
three points $P_1$, $P_{13}$ and $P_3$, and that the curve $ ( x_4, f_{n, p}(x_4) )$ for $\displaystyle  (  \frac{2 p }{n} \leq x_4 \leq  \frac{2 (p+1)}{n}   )$ lies inside the triangle given by the
three points $P_3$, $P_{23}$ and $P_2$. 
 Since the point $(u_1,  f_{n, p}(u_1) )$ is inside of one of these  triangles, we see that   the local minimum  $ f_{n, p}(u_1)  $ is positive for $n/2 \leq p \leq n-1$, 
 and thus we get our claim.
 
\medskip 
 
 {\bf Case 2.}  We obtain the equation 
 \begin{eqnarray*}  & &  \ \ 
 {x_1}{x_4} (n-p+1)+{x_3}{x_4} (n-p+1)-2 (n+1)
  {x_2}{x_4}+(p+1){x_1}{x_2}+(p+1)
  {x_2}{x_3} = 0, 
    \end{eqnarray*}

we put $x_1 =1$ and 
  we  get a following system of  equations
  \begin{eqnarray} & & 
\  -(n+p+1){x_2}^2{x_4} -(n-3 p+1) {x_3}^2{x_4}
   +(n+p+1){x_4}  +4 (n+1){x_2}{x_3}
  {x_4} \nonumber  \\
  &  & -4 (n+1){x_3}{x_4}-(p+1){x_2}
  {x_3}^2-(p+1){x_2}{x_4}^2+(p+1)
  {x_2}  = 0,  \label{5}  \\
 & & \  -(n-p+1){x_2}^2{x_4} +  (2 n-p+1){x_2}
  {x_3}^2 -  (2 n-p+1){x_2}{x_4}^2 +  (2 n-3 p-1){x_2}  \nonumber  \\
   & & +  (n-p+1) {x_3}^2{x_4}
  -  (n-p+1){x_4}-4 (n+1){x_2}{x_3}+4
   (n+1){x_2}{x_4} = 0,  \ \  \label{6}\\
 & &  \   {x_4} (n-p+1)+{x_3}{x_4} (n-p+1)-2 (n+1)
  {x_2}{x_4}+(p+1){x_2}+(p+1)
  {x_2}{x_3} = 0.  \ \  \ \  \label{7}
\end{eqnarray} 
From  equation (\ref{7}), we have 
 \begin{eqnarray} & & 
  {x_2} = \frac{({x_3}+1) {x_4} (n-p+1)}{2
   (n+1) {x_4} - (p+1) ({x_3}+1)}. \label{8}
   \end{eqnarray} 
 Now we substitute equation  (\ref{8})  into  equations (\ref{5}) and  (\ref{6}), and we obtain the following equations: 
 \begin{eqnarray}
& & F_1^{}(x_3, x_4) = -\left(3 n^3+5 n^2 p+9 n^2+2 n p^2+12 n p+10 n-2
   p^3+6 p+4\right){x_3}^2 {x_4}^2 \nonumber \\
    & & +2 \left(5 n^3+3 n^2 p+15 n^2-2 n
   p^2+4 n p+14 n+2 p^3+2 p+4\right)  {x_3} {x_4}^2   \nonumber  \\
   & & - \left(3 n^3+5 n^2 p+9 n^2+2
   n p^2+12 n p+10 n-2 p^3+6 p+4\right){x_4}^2  \nonumber  \\
   & & +2 (n+1) (p+1) (n+3 p+1) {x_3}^3 {x_4}+4 (p+1)^2  (n-p+1){x_3}^3  \nonumber  \\
   & & -2 (n+1) (p+1)  (5 n-p+5) {x_3}^2 {x_4}+4 (p+1)^2  (2 n-p+2){x_3}^2  \nonumber \\
   & & +2 (n+1) (p+1)  (n-p+1) {x_3} {x_4}^3 -2 (n+1) (p+1)  (5 n-p+5){x_3} {x_4}  \nonumber  \\
   & & +4 (p+1)^2  (n-p+1){x_3} +2 (n+1) (p+1) (n-p+1)  {x_4}^3  \nonumber  \\
   & & +2 (n+1) (p+1)  (n+3 p+1) {x_4} -2 p (p+1)^2 {x_3}^4 -2 p (p+1)^2 = 0,  \label{9} \\ 
   & &  \nonumber  \\
   & & F_2^{}(x_3, x_4) = 2 (p+1)(n-p) {x_3}^3  -2(n+1) (2 n-3 p-1){x_3}^2 {x_4}   \nonumber  \\
   & & - \left(3 n^2+4 n p+8 n-2 p^2+2 p+4\right){x_4}^2 
  -\left(3 n^2+4 n p+8 n-2 p^2+2 p+4\right){x_3} {x_4}^2  \nonumber  \\
  & &  -2 (p+1) (n+p+2) {x_3}^2 
   +4 (n+1) (2 n+p+3){x_3} {x_4}  -2 (p+1) (n+p+2) {x_3}  \nonumber  \\
   & & 
   +2 (n+1)  (2 n-p+1){x_4}^3 -2 (n+1) (2 n-3 p-1){x_4}  +2 (p+1) (n-p) = 0.  \label{10}
  \end{eqnarray} 
   We consider the resultant of the polynomials $ F_1^{}(x_3, x_4)$ and  $ F_2^{}(x_3, x_4)$ with respect to $x_4$,  which is a polynomial  of $x_3$, say $P(x_3)$. 
  We factor $P(x_3)$ as 
    \begin{eqnarray*}
& &  \ \ P(x_3) = -16 (n+1)^4 (p+1)^2 ({x_3}+1)^4 (n-p+1)^2
   (n {x_3}-n-2 p-2)\times \\
   & &  (n {x_3}-3 n+2
   p-2) ({x_3} (3 n-2 p+2)-n) ({x_3}
   (n+2 p+2)-n)\times \\
& &    \left(n^2 (3 n+4) {x_3}^4 -8 n \left(2 n^2+n p+5 n-p^2+3\right) {x_3}^3  \right. \\ 
& & 2 \left(13 n^3+8 n^2 p+36 n^2-8 n p^2+16 n
   p+40 n-16 p^2+16\right) {x_3}^2 \\
   & &  \left.  -8 n \left(2 n^2+n p+5 n-p^2+3\right) {x_3} +n^2 (3 n+4)\right). 
 \end{eqnarray*}
We denote by $Q_{n, \, p} (x_3)$ the  factor of degree $4$ in the above factorization:     
         \begin{eqnarray*}
& &  \ \ Q_{n, \, p} (x_3) =   n^2 (3 n+4) {x_3}^4 -8 n \left(2 n^2+n p+5 n-p^2+3\right) {x_3}^3  \\
& & 2 \left(13 n^3+8 n^2 p+36 n^2-8 n p^2+16 n
   p+40 n-16 p^2+16\right) {x_3}^2 \\
   & &   -8 n \left(2 n^2+n p+5 n-p^2+3\right) {x_3} +n^2 (3 n+4). 
    \end{eqnarray*}
 Now we consider  two cases:   
         
(a)\,  the case when  

\smallskip

$  (n {x_3}-n-2 p-2) (n {x_3}-3 n+2 p-2) ({x_3} (3 n-2 p+2)-n) ({x_3}
   (n+2 p+2)-n) = 0,   
$ 

\medskip    
     (b)\,  the case when $Q_{n, \, p} (x_3) = 0$.

\medskip

\noindent  We claim that we get only  K\"ahler-Einstein metrics on $ Sp(n )/ (U(p) \times U( n -p))$  in   case (a). 
  
1)  The  case when $\displaystyle x_3 = \frac{ n+ 2 p +2}{n}$. In this case  equations (\ref{9})  and  (\ref{10}) are given by 
\begin{eqnarray*} 
  & & \frac{4}{n^4} (n+p+1) (n {x_4}-2 p-2) (2 (p+1)
   (n+p+1)-n (n-p+1) {x_4})  \times\\
  & & \left((-n-p-1)
   \left(n^2+2 n-2 p^2-2 p\right)-n (n+1)
   (p+1) {x_4}\right) = 0,  \\
  & & \frac{2}{n^3} (n {x_4}-2 p-2) \left(n^2
   (n+1) (2 n-p+1){x_4}^2 -n \left(3 n^3+3 n^2 p+7 n^2+4 n p^2   \right. \right.\\
   & & \left.  \left.    +10 n p+6  n-2 p^3+2 p^2+6 p+2\right){x_4} + 2 (p+1)
   (n+p+1) \left(n^2+2 p^2+2 p\right)\right)  = 0. 
      \end{eqnarray*} 
If $  n {x_4}-2 p-2 \neq 0$, we have 
\begin{eqnarray*} 
  & & \ \  (2 (p+1) (n+p+1)-n (n-p+1) {x_4}) \times\\
  & &  \left((-n-p-1)
   \left(n^2+2 n-2 p^2-2 p\right)-n (n+1)
   (p+1) {x_4}\right) = 0,  \\ 
  & &\ \   \left(n^2
   (n+1) (2 n-p+1){x_4}^2 -n \left(3 n^3+3 n^2 p+7 n^2+4 n p^2 +10 n p+6  n-2 p^3  \right. \right.\\
   & & \left.  \left.    +2 p^2+6 p+2\right){x_4} + 2 (p+1)
   (n+p+1) \left(n^2+2 p^2+2 p\right)\right)  = 0. 
      \end{eqnarray*} 
 By taking   the resultant of these polynomials with respect to $x_4$, 
 we get 
 $$
-8 n^6 (n+1)^3 (p+1) (n-2 p) (n-p+1) (n+p+1)^4
   \left(n^2+2 n p+4 n-2 p^2+2\right),$$
  and we see that  the resultant is non-zero for $ 1 \leq p \leq n-1$ and $n \neq 2p$. Thus we get only  $\displaystyle  {x_4} = \frac{2 p  +2}{n} $ for a solution of equations (\ref{9})  and  (\ref{10}). From (\ref{8}), we see  $\displaystyle  {x_2} = \frac{ 2( n+  p+ 1)}{n}$. For $n =2 p$, we get $\displaystyle  {x_4} = \frac{ p  + 1}{p} $ and $\displaystyle  {x_4} = \frac{3 p  + 1}{p} $ as  solutions of equations (\ref{9})  and  (\ref{10}).  From (\ref{8}), we see  $\displaystyle  {x_2} = \frac{ 3 p+ 1}{p}$ and  $\displaystyle  {x_2} = \frac{  p+ 1}{p}$ respectively. 
Thus we get  K\"ahler-Einstein metrics in this case. 

2)  The  case when $\displaystyle x_3 = \frac{3 n- 2 p +2}{n}$. In this case equations (\ref{9})  and  (\ref{10}) are given by 
\begin{eqnarray*} 
 & & \ \  \frac{4}{n^4} (p+1) (n-p+1) (n {x_4}-4 n+2 p-2)
   \left(n^2 (n+1)(2 n-p+1){x_4}^2  \right. \\
   & &  +n \left(4
   n^3-11 n^2 p+n^2+10 n p^2-10 n p-4 n-2
   p^3+6 p^2-2 p-2\right) {x_4} \\
   & & \left.  -2 (p+1) (2 n-p+1) \left(3 n^2-4 n
   p+2 n+2 p^2-2 p\right) \right) = 0, \\
    & & \ \  \frac{2}{n^3} (2 n-p+1) (n {x_4}-4 n+2 p-2) \times \\
     & &  (n {x_4}-2 p-2) \left((n+1)
   n{x_4}+n^2-4 n p +2 p^2-2
   p\right)  = 0.
  \end{eqnarray*} 
If $  n {x_4}-4 n + 2 p-2 \neq 0$, we have 
\begin{eqnarray*} 
  & & \ \ n^2 (n+1)(2 n-p+1){x_4}^2 +n \left(4
   n^3-11 n^2 p+n^2+10 n p^2-10 n p-4 n-2
   p^3 \right. \\
   & & \left.  +6 p^2-2 p-2\right) {x_4}  -2 (p+1) (2 n-p+1) \left(3 n^2-4 n
   p+2 n+2 p^2-2 p\right)  = 0, \\
  & &\ \  (n {x_4}-2 p-2) \left((n+1)
   n{x_4}+n^2-4 n p +2 p^2-2
   p\right)  = 0. 
      \end{eqnarray*} 
 By taking   the resultant of these polynomials with respect to $x_4$, 
 we get 
 $$8 n^6 (n+1)^3 (p+1) (n-2 p) (n-p+1)
   \left(n^2+2 n p+4 n-2 p^2+2\right),  $$
 and we see that  the resultant is non-zero for $ 1 \leq p \leq n-1$ and $n \neq 2p$. Thus we get only  $\displaystyle  {x_4} = \frac{4 n - 2 p  + 2}{n} $ for a solution of equations (\ref{9})  and  (\ref{10}). From (\ref{8}), we see  $\displaystyle  {x_2} = \frac{ 2( n-  p+ 1)}{n}$. For $n =2 p$, we get $\displaystyle  {x_4} = \frac{ p  + 1}{p} $ and $\displaystyle  {x_4} = \frac{3 p  + 1}{p} $ as  solutions of equations (\ref{9})  and  (\ref{10}).  From (\ref{8}), we see  $\displaystyle  {x_2} = \frac{ 3 p+ 1}{p}$ and  $\displaystyle  {x_2} = \frac{  p+ 1}{p}$ respectively. 
Thus we get  K\"ahler-Einstein metrics in this case. 

3)  The  case when $\displaystyle x_3 = \frac{n}{ n+ 2 p +2}$. In this case equations (\ref{9})  and  (\ref{10}) are given by 
\begin{eqnarray*} 
 & & \ \  \frac{4 (n+p+1)}{(n+2 p+2)^4} ( (n+2 p+2) {x_4} -2 (p+1))\times  \\
 & &
(2 (p+1) (n+p+1) - (n-p+1) (n+2 p+2){x_4}) \times \\
& & 
   \left(-n^3-n^2 p-3 n^2+2 n p^2 -2 n+2 p^3+4 p^2+2 p 
   -(n+1)(p+1)(n+2 p+2){x_4}\right) = 0, \\
    & & \ \  \frac{2}{(n+2 p+2)^3}
   ( (n+2 p+2){x_4}-2 (p+1))\cdot \left(2 (p+1) (n+p+1) \left(n^2+2 p^2+2 p\right)\right.  \\
   & &  -(n+2 p+2) \left(3
   n^3+3 n^2 p+7 n^2+4 n p^2+10 n p+6 n-2
   p^3+2 p^2+6 p+2\right) {x_4} \\
   & &   \left. +(n+1)
   (2 n-p+1) (n+2 p+2)^2 {x_4}^2 \right)  = 0.
  \end{eqnarray*} 
  If $   (n+2 p+2) {x_4} -2 (p+1) \neq 0$, we have 
\begin{eqnarray*} 
& &  \ \ 
(2 (p+1) (n+p+1) - (n-p+1) (n+2 p+2){x_4}) \times \\
& & 
   \left(-n^3-n^2 p-3 n^2+2 n p^2 -2 n+2 p^3+4 p^2+2 p 
   -(n+1)(p+1)(n+2 p+2){x_4}\right) = 0, \\
& &   \ \  \left(2 (p+1) (n+p+1) \left(n^2+2 p^2+2 p\right)
   -(n+2 p+2) \left(3
   n^3+3 n^2 p+7 n^2+4 n p^2 \right.  \right.  \\
  & &   \left. \left. +10 n p+6 n-2
   p^3+2 p^2+6 p+2\right) {x_4} 
 +(n+1)
   (2 n-p+1) (n+2 p+2)^2 {x_4}^2 \right)  = 0.
  \end{eqnarray*} 
   By taking   the resultant of these polynomials with respect to $x_4$, 
 we get 
 $$-8 n^2 (n+1)^3 (p+1) (n-2 p) (n-p+1)
   (n+p+1)^4 (n+2 p+2)^4 \left(n^2+2 n p+4
   n-2 p^2+2\right),  $$
 and we see that  the resultant is non-zero for $ 1 \leq p \leq n-1$ and $n \neq 2p$. Thus we get only  $\displaystyle  {x_4} = \frac{ 2 ( p  + 1)}{n+2 p+2} $ for a solution of equations (\ref{9})  and  (\ref{10}). From (\ref{8}), we see  $\displaystyle  {x_2} = \frac{ 2( n+  p+ 1)}{n+2 p+2}$. For $n =2 p$, we get $\displaystyle  {x_4} = \frac{ p  + 1}{2 p +1} $ and $\displaystyle  {x_4} = \frac{3 p  + 1}{2 p+1} $ as  solutions of equations (\ref{9})  and  (\ref{10}).  From (\ref{8}), we see  $\displaystyle  {x_2} = \frac{ 3 p+ 1}{2 p+1}$ and  $\displaystyle  {x_2} = \frac{  p+ 1}{2 p+ 1}$ respectively. 
Thus we get  K\"ahler-Einstein metrics in this case. 

4)  The  case when $\displaystyle x_3 = \frac{n}{3 n- 2 p +2}$. In this case equations (\ref{9})  and  (\ref{10}) are given by 
\begin{eqnarray*} 
 & & \ \  \frac{4 (p+1) (n-p+1)}{(3 n-2 p+2)^4} ( (3
   n-2 p+2){x_4}-2 (2 n-p+1)) \cdot  \left(-2 (p+1) (2
   n-p+1)\times \right. \\
   & &  \left(3 n^2-4 n p+2 n+2 p^2-2
   p\right)+ (3 n-2 p+2) \left(4
   n^3-11 n^2 p+n^2  +10 n p^2-10 n p \right. \\
   & &  \left.  \left. -4 n-2
   p^3+6 p^2-2 p-2\right){x_4} +(n+1)
   (3 n-2 p+2)^2 (2 n-p+1){x_4}^2 \right)  = 0, \\
    & & \ \   \frac{2(2 n-p+1) }{(3 n-2
   p+2)^3}( (3 n-2 p+2){x_4}-2 (2
   n-p+1)) ( (3 n-2 p+2){x_4}-2 (p+1)) \times \\ 
  & &   \left((n+1)(3 n-2 p+2){x_4}+n^2 -4
   n p+2 p^2-2 p\right) = 0.
  \end{eqnarray*} 
   If $   (3 n-2 p+2){x_4}-2 (2 n-p+1) \neq 0$, we have 
\begin{eqnarray*} 
& &  \ \ 
(2 (p+1) (n+p+1) - (n-p+1) (n+2 p+2){x_4}) \times \\
& & 
   \left(-n^3-n^2 p-3 n^2+2 n p^2 -2 n+2 p^3+4 p^2+2 p 
   -(n+1)(p+1)(n+2 p+2){x_4}\right) = 0, \\
& &   \ \  \left(2 (p+1) (n+p+1) \left(n^2+2 p^2+2 p\right)
   -(n+2 p+2) \left(3
   n^3+3 n^2 p+7 n^2+4 n p^2 \right.  \right.  \\
  & &   \left. \left. +10 n p+6 n-2
   p^3+2 p^2+6 p+2\right) {x_4} 
 +(n+1)
   (2 n-p+1) (n+2 p+2)^2 {x_4}^2 \right)  = 0.
  \end{eqnarray*} 
   By taking   the resultant of these polynomials with respect to $x_4$, 
 we get 
$$ 8 n^2 (n+1)^3 (p+1) (n-2 p) (3 n-2 p+2)^4
   (n-p+1) \left(n^2+2 n p+4 n-2 p^2+2\right),$$ 
 and we see that  the resultant is non-zero for $ 1 \leq p \leq n-1$ and $n \neq 2p$. Thus we get only  $\displaystyle  {x_4} = \frac{2 (2 n-p+1)}{3 n-2 p+2} $ for a solution of equations (\ref{9})  and  (\ref{10}). From (\ref{8}), we see  $\displaystyle  {x_2} = \frac{ 2( n -  p+ 1)}{3 n-2 p+2}$. For $n =2 p$, we get $\displaystyle  {x_4} = \frac{ p  + 1}{2 p +1} $ and $\displaystyle  {x_4} = \frac{3 p  + 1}{2 p+1} $ as  solutions of equations (\ref{9})  and  (\ref{10}).  From (\ref{8}), we see  $\displaystyle  {x_2} = \frac{ 3 p+ 1}{2 p+1}$ and  $\displaystyle  {x_2} = \frac{  p+ 1}{2 p+ 1}$ respectively. 
Thus we get  K\"ahler-Einstein metrics in this case.

\medskip

Now we consider the case (b), that is, the case when  $Q_{n, \, p} (x_3) = 0$. 
We compute a Gr\"obner basis of $\{ F_1^{}(x_3, x_4), \ F_2^{}(x_3, x_4), \ Q_{n, \, p} (x_3) \}$  using the lex order with $ x_4 >  x_3$.  We can find  the following polynomials in the Gr\"obner basis: 
  \begin{eqnarray} 
& &  \ \ \   Q_{n, \, p} (x_3),   \nonumber \\
& & \ \ \    n^2 (3 n+4){x_3}^3 -n \left(19 n^2+8 n p+44 n-8
   p^2+24\right){x_3}^2  \nonumber \\
   & & + 3\left(29
   n^3+32 n^2 p+92 n^2-24 n p^2+40 n
   p+96 n-32 p^2+2\right){x_3}  \nonumber  \\
   & & -n \left(13 n^2-8 n
   p+20 n-8 p+8\right)-8 n (n+1) (2 n-p+1) {x_4},  \label{11} \\
   & &  \ \ \  n {x_3}^2 -2 (n+2 p+2){x_3}+ n ({x_3}+1) {x_4}.  \label{12} 
      \end{eqnarray}  
      
  From (\ref{12}) we have 
  \begin{eqnarray} 
  x_4 = \frac{- n {x_3}^2 + 2 {x_3} (n+2 p+2)  - n}{n ({x_3}+1)}.  \label{13} 
  \end{eqnarray}  
       We substitute equation  (\ref{13})  into  equation (\ref{8}), and we obtain   
        \begin{eqnarray}  
  x_2 =  \frac {({x_ 3} + 1) (-n + p - 1) \left ( - 
     n {x_ 3}^2 + 2 {x_ 3} (n + 2 p + 2)  - 
     n \right)} {  n (2 n + p + 3) {x_ 3}^2 -2  \left (2 n^2 + 3 n p + 5 n + 4 p + 
      4 \right){x_ 3} + n (2 n + p + 3)}.  \label{14}   
  \end{eqnarray}  
   
Let $a_k$ ($ k = 0, \cdots, 4$) denote the coefficients of the polynomial  $ Q_{n, \ p} (x_3) $. Then we have that $a_k = a_{4-k}$.  Thus we see that if  the equation $ Q_{n, \ p} (x_3) = 0$ has a solution $ x_3  = \alpha$, then so is $\displaystyle  x_3  = \frac{1}{\alpha}$. 

We claim that the equation $ Q_{n, p} (x_3) = 0$ has four different positive solutions. It is enough to see that $ Q_{n,  p} (x_3) = 0$ has two different solutions
between $ 0 < x_3  < 1$. 

Note that $ Q_{n,  p} (0) = n^2 (3 n+4) > 0 $ and  $ Q_{n, p} (1) = 32 (p+1) (n-p+1) > 0 $. 
 \begin{eqnarray*}
  & & 
  Q_{n, p} \left(\frac{1}{2}\right) =
\frac{1}{16} \left(-5 n^3-16 n^2 p-44 n^2+16 n
   p^2+128 n p+80 n-128 p^2+128\right) \\
& &   = -5(n-p-1)^3+(-31 p-59) (n-p-1)^2+\left(-31 p^2-22 p-23\right) (n-p-1) \\
& &    -5 p^3-75
   p^2+89 p+159
   \\ 
   & & 
   = -5(n-p-1)^3+(-31 p-59) (n-p-1)^2+\left(-31 p^2-22 p-23\right) (n-p-1) \\
 & &   -5 (p-2)^3-105 (p-2)^2-271 (p-2)-3 < 0  \quad \mbox{ for } \ \  2 \leq p \leq n-1. 
       \end{eqnarray*}  
       For $ p =1$, we have 
        \begin{eqnarray*}
  & & 
        Q_{n,  1} \left(\frac{1}{2}\right) = -\frac{1}{16} n \left(5 n^2+60 n-224\right)
         = -\frac{1}{16} \left(5 (n-3)^2+90 (n-3)+1\right) n < 0 
      \end{eqnarray*} 
      for $ n \geq 3$.  Thus we get our claim. 

Now, for a solution of the equation $ Q_{n,  p} (x_3) = 0$, we get  values of $x_4$ and  $x_2$ from (\ref{13}) and  (\ref{14}), and thus  we get four triples of solutions $\{x_2^0, x_3^0, x_4^0 \}$ of 
the system of equations  (\ref{5}),   (\ref{6}),   (\ref{7}) in case (b). 

We claim that only two  triples of the solutions  of 
the system of equations  (\ref{5}),   (\ref{6}),   (\ref{7}) have the property that $x_2^0 > 0$, $x_3^0 > 0$, $x_4^0 > 0$ in case (b). 

Consider the resultant of the polynomials 
$$Q_{n, p} (x_3) \ \ \mbox{and}  \ \  n {x_3}^2 -2 (n+2 p+2){x_3}+ n ({x_3}+1) {x_4} $$
   with respect to $x_3$,  which is a polynomial  of $x_4$, say  $S_{n, \, p} (x_4)$.
 Then we get that, up to a scalar multiple $32 n^4 ( n+p+1)$,  
 \begin{eqnarray*}
& &  \ \ S_{n, \, p} (x_4) =  n^2 (n+1) (3 n+4) (2 n-p+1) {x_4}^4  \\
   & & + 4 n (n+1) (2n-p+1) \left(n^2-4 n p-2 p^2-8 p-2\right) {x_4}^3  \\
   & & +2  \left(n^5-19
   n^4 p-11 n^4+36 n^3 p^2-18 n^3 p-30 n^3+22 n^2
   p^3+130 n^2 p^2+54 n^2 p \right. \\
   & & \left.  -22 n^2-16 n p^4+4 n
   p^3+108 n p^2+68 n p-4 n-16 p^4-16 p^3+16 p^2+16
   p\right){x_4}^2  \\
   & & -8 (p+1)  (n-2 p) (n+p+1) \left(n^2-6 n
   p-2 n+2 p^2-4 p-2\right){x_4} \\
   & & +8 (p+1)^2 (n-2 p)^2 (n+p+1).    
         \end{eqnarray*} 
         From (\ref{14}) we get a polynomial 
            \begin{eqnarray*}  
  & & \ \ P_{n, p}(x_2, x_3) =   x_2 \cdot ( n (2 n + p + 3) {x_ 3}^2 -2  \left (2 n^2 + 3 n p + 5 n + 4 p + 
      4 \right){x_ 3} + n (2 n + p + 3)) \\
  & &     - ({x_ 3} + 1) (-n + p - 1) \left ( - 
     n {x_ 3}^2 + 2 {x_ 3} (n + 2 p + 2)  - 
     n \right).  
  \end{eqnarray*}  
We also consider the resultant of the polynomials $Q_{n, p} (x_3)$ and $P_{n, p}$  with respect to $x_3$, which is a polynomial  of $x_2$, say  $T_{n, \, p} (x_2)$.
Then we get that, up to a scalar multiple $1024 n^4 (n+1) (n-p+1)^2 (p+1)^2 ( n+p+1)^3$,  
 \begin{eqnarray*} 
 & &  \ \ T_{n,  p} (x_2) =         
 n^2 (n+1) (3 n+4) (n+p+1) {x_2}^4  \\
   & & -4 n (n+1)  (n+p+1) \left(5
   n^2-8 n p+8 n+2 p^2-8 p+2\right){x_2}^3 \\
      & & +2 \left(24 n^5-55 n^4 p+89 n^4+6
   n^3 p^2-190 n^3 p+116 n^3+42 n^2 p^3+46 n^2 p^2-222 n^2 p \right.  \\
    & & \left.  +62 n^2-16
   n p^4+60 n p^3+60 n p^2-100 n p+12 n-16 p^4+16 p^3+16 p^2-16
   p\right)  {x_2}^2 \\
   & & -8 (n-2 p) (n-p+1) (2
   n-p+1) \left(3 n^2-2 n p+6 n-2 p^2-4 p+2\right) {x_2} \\
   & & +8 (n-2 p)^2 (n-p+1)^2 (2 n-p+1).  
          \end{eqnarray*}  
Note that   $T_{n,  p} (x_2) = S_{n, \, n - p} (x_2)$ by a direct computation. 

\medskip

Now we claim that 

 \smallskip
 
(I) \  for $\displaystyle 1 \leq p  <  \frac{n}{2}$, the equation $S_{n, \, p} (x_4) = 0$ has  
 two    different positive solutions and  two   different negative solutions
 
  \smallskip
      
      and 
   
    \smallskip
    
(II)   \  for $\displaystyle   \frac{n}{2} < p \leq n-1$, the equation $S_{n, \, p} (x_4) = 0$ has   four   different positive solutions. 

\medskip

Note that,  for $ x_4 = 0$,  we have $S_{n, \, p} (0) = 8 (p+1)^2 (n-2 p)^2 (n+p+1) > 0$ for $\displaystyle  p \neq  \frac{n}{2}$,  

\noindent   and for $\displaystyle    x_4 = \frac{2 (p+1)}{n}$,  we have 
     \begin{eqnarray*}
& &  \ \ \displaystyle S_{n, \, p} \left( \frac{2 (p+1)}{n} \right) = -\frac{1}{n^2}16 (p+1)^2 (n-p+1)^2 \left(n ( p-1) +4 p^2+4
   p\right) < 0. 
     \end{eqnarray*}
   For $\displaystyle    x_4 = \frac{2 (2 p-n)}{n}$,  we have 
     \begin{eqnarray*}
     \displaystyle S_{n, \, p} \left( \frac{2 (2 p- n)}{n} \right) =\frac{8}{n} (n-2 p)^2 (n-p+1)^2 \left(5 n^2-9 n p+3
   n+4 p^2-4 p\right). 
    \end{eqnarray*} 
Note that    
 \begin{eqnarray*}
  5 n^2-9 n p+3 n+4 p^2-4 p = 5 (n-p)^2+(p+3) (n-p)-p \\
   > (n-p)^2+(p+3)\cdot 1 -p =(n-p)^2 + 3 > 0. 
    \end{eqnarray*} 
    Thus  we get  $
     \displaystyle S_{n, \, p} \left( \frac{2 (2 p- n)}{n} \right) > 0.  $

\noindent  Now  for $\displaystyle    x_4 = \frac{ (2 p-n)}{n}$,  we have 
     \begin{eqnarray*}
& &  \ \ S_{n, \, p} \left( \frac{2 p- n}{n} \right) = -\frac{(n-2 p)^2}{n^2} \left(5 n^4 p+9 n^4-4 n^3 p^2+3
   n^3 p+33 n^3+4 n^2 p^2  \right.\\
   & &  \left. -20 n^2 p+40 n^2-20 n
   p^3+4 n p^2-32 n p+16 n+16 p^4-16
   p \right) \\
   & & =  -\frac{(n-2 p)^2}{n^2} \left( \left(16 p^2+39 p+33\right) (n-p)^3+\left(18
   p^3+67 p^2+79 p+40\right) (n-p)^2 +\left(8 p^4 \right. \right. \\
   & & \left. \left. +33 p^3+63 p^2+48 p+16\right) (n-p)+(5
   p+9) (n-p)^4+p^5+12 p^4+17 p^3+8 p^2 \right)  < 0.  
    \end{eqnarray*}

For  $\displaystyle 1 \leq p  <  \frac{n}{2}$,  we have $2 p - n  < 0  $  and hence 
$2(2p -n) < 2p -n < 0 < 2(p+1)$.         
Thus  we see that the equation $S_{n, \, p}(x_4) = 0$ has the four solutions $x_4^1$, $x_4^2$, $x_4^3$, $x_4^4$ with 

$$ \displaystyle  \frac{2(2 p- n)}{n} < x_4^1 <  \frac{2 p- n}{n} < 
 x_4^2 <  0 < x_4^3 <   \frac{ 2(p+1)}{n} < x_4^4. $$  
 
 For  $\displaystyle   \frac{n}{2} < p \leq n-1$,   we have $2 p - n  > 0  $. 
 Since $  2(2 p- n) - 2(p + 1) = 2(p -n) < 0$, we have $0 <  2p -n < 2(2p -n) <  2(p+1)$. 
 Thus  we see that the equation $S_{n, \, p}(x_4) = 0$ has the four solutions $x_4^1$, $x_4^2$, $x_4^3$, $x_4^4$ with 
$$ \displaystyle 0 <  x_4^1 <  \frac{2 p- n}{n} < x_4^2 <  \frac{2(2 p- n)}{n} < 
 x_4^3 <    \frac{ 2(p+1)}{n} < x_4^4. $$  

Noting that   $T_{n, \, p} (x_2) = S_{n, \, n - p} (x_2)$, we also get that 

(III) \  for $\displaystyle 1 \leq p  <  \frac{n}{2}$, the equation $T_{n, \, p} (x_2) = 0$ has   four   different positive solutions. 
 
  \smallskip
      
      and 
   
    \smallskip
    
(IV)   \  for $\displaystyle   \frac{n}{2} < p \leq n-1$, the equation $T_{n, \, p} (x_2) = 0$ has  
 two   different positive solutions and  $2$   different negative solutions. 

Combining the statements (I), (II), (III), and (IV) we 
get exactly  two non-K\"ahler Einstein metrics on $Sp(n)/ (U(p) \times U( n -p))$ and this completes the proof.

\end{document}